\documentstyle[amssymb]{amsart}
\newcommand{\REF}[1]{[#1]}
\newcommand{\noproof}{{\unskip\nobreak\hfill\penalty50\hskip2em
    \hbox{}\nobreak\hfill$\blacksquare$\parfillskip=0pt
    \finalhyphendemerits=0}}
\newcommand{\stopproof}{\noproof\bigskip}
\newcommand{\ged}{\ge^*}
\newcommand{\led}{\le^*}

\newcommand{\proof}{{\bf Proof:} \ }

\newcommand{\lcard}{\, \mid \!}
\newcommand{\rcard}{\! \mid \,}

\newcommand{\card}[1]{\lcard #1 \rcard}
\newtheorem{theor}{Theorem}[section]
\newtheorem{propo}{Proposition}[section]

\newtheorem{lemma}{Lemma}[section]

\author{Reinhard Diestel}
\address{Fakultat fur Mathematik, Universitat Bielefeld \\
SFB 343 \\
D4800 Bielefeld 1 \\ Germany }
\author{Saharon Shelah}
\address{Institute of Mathematics\\Hebrew University\\Jerusalem,
Givat Ram,  Israel}
\author{Juris Stepr\={a}ns}
\address{Department of Mathematics, York University \\
4700 Keele Street \\
North York, Ontario \\ Canada  M3J 1P3}
\title{Characterizing Dominating Graphs}
\thanks{The second author's research is supported by the Israel Research Fund
while the third author is supported by NSERC. This is publication 443 in the
second author's catalogue of publications}
\begin{document}
    \maketitle
\begin{abstract}
    A graph is called dominating if its vertices can be labelled with
integers in such a way that for every function $f\colon\,\omega\to
\omega$ the graph contains a ray whose sequence of labels eventually
exceeds~$f$. We obtain a characterization of these graphs by producing
a small family of dominating graphs with the property that every
dominating graph must contain some member of the family.
 \end{abstract}

\section{Introduction}
If $f$ and $g$ are functions from $\omega$ to $\omega$, we write
$f\ged  g$ and say that $f$ {\em dominates}~$g$ if the set $\{n\in
\omega : f(n) < g(n)\}$ is finite. A family $\cal{F}$ of functions from
$\omega$ to $\omega$ is called a {\em dominating family} if every
function $g\colon\,\omega \rightarrow \omega$ is dominated by some
$f\in\cal{F}$. The least cardinality of a dominating family is denoted
by~$\frak{d}$.

Similarly, a family $\cal{F}$ of functions from $\omega$ to $\omega$ is
called {\em bounded} if there exists a function $g\colon\,\omega
\rightarrow \omega$ which dominates every $f\in\cal{F}$; if no such
function exists, $\cal F$ is {\em unbounded}.  The least cardinality
of an unbounded family is denoted by~$\frak{b}$.

It is well known and easy to show that $\omega < \frak b \le \frak d
\le 2^\omega$. Depending on the axioms of set theory assumed, $\frak
b$ and $\frak d$ may or may not coincide, and it is consistent that
both are less than~$2^\omega$. Properties of these and related
cardinals have been studied widely in the literature; see the article
by Vaughan in~\REF3.

Taking a different approach to considering merely the cardinalities of
boun\-ded families of functions, Halin (see~\REF2) introduced the
notion of a bounded graph: a graph is called {\em bounded} if for every
labelling of its vertices with integers the labellings along its
rays---its one-way infinite paths---form a bounded family. Thus, the
family of functions considered is constrained not by cardinality but
by imposing an intersection pattern on its members. A long-standing
conjecture of Halin, known as the `bounded graph conjecture', said
that the bounded graphs are characterized by the exclusion of four
simple types of unbounded graph; this conjecture was recently proved
in~\REF1.

In this paper we prove an analogous result for dominating graphs; a
graph is called {\em dominating} if its vertices can be labelled with
integers in such a way that the labellings along its rays form a
dominating family of functions. We show that a graph is dominating if
and only if it contains one of three specified prototypes of a
dominating graph.

\medbreak

As usual, a {\em graph} will be thought of as a symmetrical binary
relation on some underlying set, its set of {\em vertices}. Thus, a
graph on a set $X$ is a subset of the set $[X]^2$ of unordered pairs
of~$X$, called its {\em edges}. Two graphs will be called {\em
disjoint} if and only if their vertex sets are disjoint. If $G$ is a
graph on~$X$, then $G'$ is a {\em subgraph} of $G$ if $G'$ is a subset
of $G\cap [X']^2$ for some $X'\subset X$.

The {\em degree} of a vertex is the number of edges containing it. If
$m\in\omega+1$, a {\em path} of {\em length}~$m$ in a graph $G$ on~$X$
is a one-to-one function $P\colon\,m \rightarrow X$ such that
$\{P(n-1),P(n)\}\in G$ whenever $0 < n < m$. Often, the image of a path
will be confused with the path itself; for example, a vertex $x$ will
be said to be `on'~$P$ when what is really meant is that $P(n) = x$
for some $n\in m$. With some abuse of notation we shall say that
$\{P_i : i\in I\}$ is a family of {\em disjoint paths from} (or: {\em
starting at})~$x$ if $P_i (0) = x$ for every~$i$ and no vertex other
than $x$ is on both $P_i$ and $P_j$ if $i\ne j$. Similarly we may
speak of a family of `disjoint' paths ending at~$x$, or of a family of
`disjoint' paths from $x$ to~$y$ when $x$ and $y$ are two fixed
vertices.

A path of infinite length will be called a {\em ray}. Thus, more
formally, a graph on $X$ is dominating if and only if there exists a
labelli
ng $L\colon\,X\to \omega$ such that for every $f\colon\, \omega
\to \omega$ there is a ray $R\colon\,\omega\to X$ with $f\led L\circ
R$.

A graph in which any two vertices can be connected by a unique path is
a {\em tree}. The tree in which every vertex has countably infinite
degree is denoted by~$T_\omega$. A tree~$T$ is called a {\em
subdivision} of~$T_\omega$ if each vertex of~$T$ has either degree~2 or
countably infinite degree, and every ray in $T$ contains a vertex that
has infinite degree in~$T$. The vertices of infinite degree in $T$ are
its {\em branch vertices}, the vertices of degree~2 its {\em
subdividing vertices}.

If  $T$ is a subdivision of~$T_\omega$, there is a natural bijection
$\phi$ from the vertices of $T_{\omega}$ to the branch vertices of~$T$
such that if $x,y$ form an edge of~$T_{\omega}$ then the unique path
in $T$ joining $\phi(x)$ to $\phi(y)$ contains no other branch vertex
of~$T$; identifying the vertices of $T_{\omega}$ with their images
under~$\phi$, we may call such a path in $T$ a {\em subdivided edge}
(of~$T_{\omega}$) {\em at}~$\phi(x)$.

A subdivision $T$ of~$T_\omega$ will be called {\em uniform} if it has
a branch vertex~$r$, called its {\em root}, such that whenever $x$ is a
branch vertex, all the subdivided edges at $x$ that are {\em not}
contained in the unique path from $x$ to~$r$ have the same length.

It is not difficult to see\REF{1, \S4} that the edges of a $T_\omega$
may be enumerated in such a way that, for every edge other than the
first edge, one of its two vertices also belongs to an edge preceding
it in the enumeration. Such an enumeration will be called a {\em
standard construction} of~$T_\omega$. As a typical (if trivial)
application of this tool, consider the task of constructing a
$T_\omega$ subgraph in some given graph every vertex of which has
infinite degree: at each step, we will have specified only a finite
portion of our $T_\omega$, so we will always be able to add the next
edge as required.

\section{Examples of dominating graphs}
In this section we look at some typical dominating graphs, including
those needed to state our characterization theorem.

Since supergraphs of dominating graphs are again dominating, our aim
will be to find dominating graphs which are minimal, in the sense that
any subgraph that does not itself contain a copy of the original graph
is no longer dominating. A trivial example of such a minimal dominating
graph is given by any graph that is the union of $\frak d$ disjoint
rays:

\begin{propo}
    If a graph \label{rays}is the union of $\frak{d}$ disjoint rays
then it is dominating.
 \end{propo}
 \proof
 Label each ray by a different member of some dominating family of
functions.
 \stopproof

So how about countable graphs? Clearly, a complete graph (one in which
every pair of vertices is an edge) on a countably infinite set is
dominating: just label its vertices injectively. In the same way we see
that a~$T_\omega$ (which is `smaller' than a complete infinite graph)
is dominating.

An arbitrary subdivision of~$T_\omega$ is not necessarily dominating.
Indeed, consider any enumeration $e\colon\, \omega\to T_\omega$ of the
edges of~$T_\omega$. For each~$n\in\omega$ subdivide $e(n)$ exactly
$n$ times, so that the resulting subdivided edge is a path of length
$n+2$. Call this tree~$T$. To see that $T$ is not dominating, let $L$
be any labelling of its vertices. Let $H\colon\,\omega\to\omega$ be any
increasing function satisfying $H(n) > \max\,\{L(x) : x\in e(n)\}$ for
all~$n\in\omega$. We show that, for any ray~$R$ in~$T$ and any
$i\in\omega$, there exists a $k > i$ such that $H(k) > L(R(k))$ (so
$H$ is not dominated by $L\circ R$). Given such $R$ and~$i$, choose
$j,k\in \omega$ with $i < j < k$ so that $\{R(j), R(k)\} = e(n)$ for
some~$n$, and so that $U = \{R(\ell) : j\le\ell\le k\}$ contains no
other branch vertex of~$T$. Then $R$ traces out the subdivided
edge~$e(n)$, and in particular we have $k\ge \card U = n + 2$. Since
$H$ is increasing and $H(n) > L(R(k
))$ by definition of~$H$, this
gives $H(k)\ge H(n) > L(R(k))$ as desired.

Uniform subdivisions of~$T_\omega$, on the other hand, are easily seen
to be dominating:

\begin{propo}
    Uniform subdivisions of~$T_\omega$ are dominating\label{type1}.
\end{propo}
\proof
Let $T$ be a uniform subdivision of~$T_\omega$, with vertex set~$X$
and root~$r$. Let $L\colon\,X\to\omega$ be any injective labelling; we
show that for every function $f\colon\omega\to\omega$ there is a ray
$R\colon\omega\to X$ such that $f\led  L\circ R$.

We define $R$ inductively, choosing its subdivided edges one at a time.
(Recall that any ray in a subdivision of $T_\omega$ contains infinitely
many branch vertices, and is thus a concatenation of paths that are
subdivided edges of the~$T_\omega$.) Let $R(0) = r$. Suppose now that
$R(n)$ has been defined for every $n\le m$, and that $R(m)$ is a
branch vertex. Then all the (infinitely many) subdivided edges
at~$R(m)$ that are not contained in the portion of $R$ defined so far
have the same length~$\ell$, and so we can find one of them, $P$~say,
such that $L(P(i)) \ge f(m+i)$ whenever $0 < i < \ell$. Setting
$R(m+i) = P(i)$ for these~$i$, we see that $L(R(m+i))\ge f(m+i)$;
moreover, $R(m+\ell-1)$ is again a branch vertex of~$T$. This
completes the induction step, and hence the construction of~$R$. Since
$L(R(n))\ge f(n)$ for every $n > 0$, we have $f\led L\circ R$ as
required.
 \stopproof

How many disjoint copies of {\em arbitrary} subdivisions of~$T_\omega$
are needed to make a dominating graph? By Proposition~\ref{rays},
$\frak d$~copies will certainly do, since each of them contains a ray.
Our next proposition says that, in fact, $\frak b$~copies suffice.

\begin{propo}
    If a graph is the \label{type2}union of $\frak{b}$ disjoint
subdivisions of~$T_\omega$, then it is dominating.
\end{propo}
\proof
 Let $\{f_{\xi} : \xi\in\frak{b}\}$ be an unbounded family of
increasing functions from $\omega$ to $\omega$. Let $\{G_{\xi}:
\xi\in\frak{b}\}$ be a family of $\frak{b}$ disjoint subdivisions
of~$T_\omega$, and let $G_{\xi}$ have vertex set $X_{\xi}$ and
root~$r_\xi$. We show that $G = \bigcup\{G_{\xi} : \xi\in
\frak{b}\}$ is dominating.

For each branch vertex $x$ of $G_{\xi}$, let $N(x)$ be the set
of all branch vertices $y$ that are not contained in the unique
path from $r_{\xi}$ to~$x$ and which are joined to $x$ by a
subdivided edge (i.e.\ by a path not containing any other branch
vertices). Let $S(x)$ denote the union of these $x$--$y$ paths;
thus, $S(x)$~consists of all the paths from $x$ to a vertex
in~$N(x)$. For $y\in N(x)$ we denote the length of the path from
$r_{\xi}$ to $y$ by~$K(y)$.

Let us define a labelling $L$ on $G$ to witness that $G$ is
dominating. For each~$\xi$, we fix $L(r_{\xi})$ arbitrarily,
and then define $L$ separately on each set $S(x)\setminus\{x\}$
for all the other branch vertices~$x$ of~$G_\xi$. There are two
cases to consider. If infinitely many $y\in N(x)$ have the same
value of $K(y)$, we let $L\upharpoonright (S(x)\setminus\{x\})$
be an arbitrary one-to-one mapping. Otherwise, we choose for each
$y\in N(x)$ some $y^+\in N(x)$ such that $K(y) < K(y^+)$; then,
for each $z\ne x$ on the path from $x$ to~$y$, we set $L(z) =
f_{\xi}(K(y^+)$.

To show that $G$ is dominating, let $f\colon\, \omega\to \omega$
be given, without loss of generality increasing. As in the proof
of Proposition \ref{type1} we inductively define a ray $R\colon\,
\omega\to X_{\xi}$ starting at~$r_\xi$, so that $f\led L\circ
R$; here $\xi$ is chosen so that $f_{\xi}\not\led f$. Since $f$
is increasing, it suffices to show that for every branch vertex
$x$ of $G_\xi$ there is some $y\in N(x)$ such that $L(z) \ge
f(K(y))$ for each $z\ne x$ on the path from $x$ to~$y$; we may
then choose the path from $x$ to $y$ as the next segment
for~$R$.

If infinitely many $y\in N(x)$ have the same value of~$K(y)$,
say~$k$, then $L$ is injective on $S(x)\setminus \{x\}$; since
$f$ takes only finitely many values on the first $k+1$ integers
,
we can easily find $y$ as desired. If not, then each $y\in N(x)$
has been assigned some $y^+\in N(x)$. Pick $y'\in N(x)$, find an
$i \ge K(y')$ such that $f_\xi (i) > f(i)$, and choose $y\in
N(x)$ with maximal $K(y) \le i$. Then $K(y) \leq i < K(y^+)$.
For each $z\ne x$ on the path from $x$ to $y$ we have
 $$L(z) = f_{\xi}(K(y^+)) \geq f_{\xi}(i) \geq f(i) \geq f(K(y))$$
as desired.
 \stopproof

\section{A characterization of dominating graphs}
We now come to prove our main result, the following characterization
of dominating graphs.
\begin{theor}
    A graph $G$ is dominating if and only if it satisfies one of the
following three conditions:\label{main}
\begin{enumerate}
    \item $G$ contains a uniform subdivision of~$T_\omega$;
    \item $G$ contains $\frak{b}$ disjoint subdivisions of~$T_\omega$;
    \item $G$ contains $\frak{d}$ disjoint rays.
\end{enumerate}
\end{theor}

Note that if $\frak b = \frak d$ then (2) above is redundant, since
$\frak d$ disjoint subdivisions of~$T_\omega$ contain $\frak{d}$
disjoint rays.

\medskip

The bulk of the proof of Theorem~\ref{main} is divided up into several
lemmas. We shall consider these lemmas in turn, and then complete
the formal proof of the theorem.

Our first lemma is an easy consequence of the fact that there is no
infinite decreasing sequence of ordinals; its proof is left to the
reader.

\begin{lemma}
    If $\rho$ is an ordinal-valued function on~$\omega$, then there
exists some $n_0\in\omega$ such that for every $n \ge n_0$ there is an
$m > n$ with $\rho(m) \geq \rho(n)$.\label{tail}\noproof
 \end{lemma}

The next three lemmas make up most of the proof of Theorem~\ref{main}.

\begin{lemma}
    If $\card{X} < \frak{b}$, then any dominating graph on $X$
    \label{t1sg} contains a uniform subdivision of~$T_\omega$.
\end{lemma}
\proof
Let $G$ be a graph on~$X$, where $\card{X} < \frak{b}$. The basic idea
of the proof is recursively to define a rank function $\rho$ on some
or all of the vertices of~$G$, with the following property. If any
vertex remains unranked, i.e.\ i f the recursion ends before $\rho$ is
defined on all of~$X$, then $G$ contains a uniform subdivision
of~$T_\omega$; if $\rho$ gets defined for every vertex, then $G$ is not
dominating.

For the definition of~$\rho$, we first define subsets $\Sigma_\xi$
of~$X$, as follows. Let $\Sigma_0$ be the set of vertices $x\in X$
that have finite degree in~$G$. For $\xi > 0$, let $\Sigma_{\xi}$ be
the set of vertices $x\in X$ such that, for every $m\in \omega$, any
family of disjoint paths of length~$m$ starting at $x$ and ending at a
vertex $y\not\in \bigcup_{\zeta\in\xi} \Sigma_{\zeta}$, is finite.
Note that if $\zeta < \xi$, then $\Sigma_{\zeta}\subset \Sigma_{\xi}$.
Finally, for each $x\in X$, define $\rho(x)$ to be the least $\xi$
such that $x\in\Sigma_{\xi}$; if no such $\xi$ exists, let $\rho(x)$
remain undefined.

It is not difficult to see that if there is some $x\in X$ such that
$\rho(x)$ is not defined then $G$ contains a  uniform subdivision
of~$T_\omega$. Indeed, if $\rho(x)$ has remained undefined then, by
definition of~$\rho$, there exists an infinite set of disjoint paths
from~$x$ in~$G$, all of the same length, and ending in vertices for
which $\rho$ is also undefined. Following the standard construction
of~$T_\omega$, it is easy to build a uniform subdivision of $T_\omega$
from all these paths: at each point of the construction, only finitely
many vertices have been used, but there is an infinite set of disjoint
paths from which the next subdivided edge can be chosen.

Let us assume from now on that $\rho(x)$ is defined for all $x\in X$,
and show that $G$ is not dominating. Let $L\colon\,X\rightarrow
\omega$ be any labelling. Assuming the Claim below (which will be
proved later), we shall find a function $H\colon\,\omega \to\omega$
which is not dominated by $L\circ R$ for any ray $R$ in~$G$.

Let a path $P$ from $x$ to $y$ in $G$ be called {\em upward} if
$\rho(y) = \max\,\{\rho(z) : z\in P\}$.

\proclaim Claim.
{\it For each $x \in X$ and $m\in\omega$, there are only finitely many
vertices $y\in X$ such that $G$ contains an upward path of length $m+1$
from $x$ to~$y$.}

>From the claim it follows that we may define, for each $x\in X$, a
function $Q_x\colon\,\omega\rightarrow \omega$ such that $Q_x(m) >
L(y)$ for any $m\in\omega$ and any vertex $y$ to which $x$ can be
linked by an upward path of length~$m+1$.  By our hypothesis that
$\card{X} < \frak{b}$, there exists a function
$H\colon\,\omega\rightarrow\omega$ which dominates each of the
functions~$Q_x$. Redefining $H(n)$ as $\max\,\{H(k) : k\leq n\}$ if
necessary, we may assume that $H$ is increasing.

Now let $R$ be any ray in~$G$; it suffices  to show that $H \not\leq^*
L\circ R$. By Lemma~\ref{tail}, we may find an infinite increasing
sequence $\{k_i : i\in \omega\}$ such that $\rho(R(k_i))\leq
\rho(R(k_{i+1}))$ for each~$i$, and $\rho(R(j)) < \rho(R(k_i))$
whenever $k_i < j < k_{i+1}$. Note in particular that, for each~$i$,
the part of $R$ that connects $R(k_0)$ with $R(k_i)$ is an upward path
of length $k_i - k_0 + 1$.

Since $H$ dominates~$Q_{R(k_0)}$, there is some $K\in \omega$ such
that $Q_{R(k_0)}(k) \leq H(k)$ for all  $k\geq K$. But then
 $$L(R(k_i)) < Q_{R(k_0)}(k_i - k_0) \leq H(k_i - k_0) \leq H(k_i)$$
for all $i$ with $k_i - k_0 \geq K$, by definition of~$Q_{R(k_0)}$.
Thus $L\circ R$ fails to dominate~$H$, as required.

\medskip

Hence all that remains to be proved is the Claim. Suppose the contrary,
and consider a vertex~$x$, an integer~$m$, and an infinite set $\{y_n :
n\in \omega\}$ such that for each $n$ there is an upward path $P_n$ of
length $m+1$ from $x$ to~$y_n$. Choose $k \le m$ maximal so that there
exist a vertex $z$ and an infinite set ${\cal P}\subset \{P_n : n\in
\omega\}$ such that $P(k) = z$ for every $P\in{\cal P}$. (Note that $k$
exists, because every~$P_n$ starts in~$x$.) We now select an infinite
sequence $\{P_{n_i} : i\in\omega\}$ of paths in~${\cal P}$ so that
any two of these are disjoint after~$z$; since each $P_n$ is an upward
path, and hence $\rho(z) \leq \rho(y_n)$ for every~$n$, this will
contradict the definition of~$\rho$.

Let $P_{n_0}$ be any path from~$\cal P$. Now suppose $P_{n_0}, \ldots,
P_{n_i}$ have been chosen, and let $U$ be the union of their vertex
sets. By the maximality of~$k$, there are at most finitely many paths
in ${\cal P}$ that contain a vertex from~$U$ after~$z$; let
$P_{n_{i+1}}$ be any other path from~$\cal P$. It is then clear that
the full sequence $\{P_{n_i} : i\in\omega\}$ has the required
disjointness property.
 \stopproof

\begin{lemma}
    If $\card{X} < \frak{d}$, then any dominating graph on $X$
    \label{t2sg} contains a subdivision of~$T_\omega$.
\end{lemma}
\proof
Let $G$ be a graph on~$X$, where $\card{X} < \frak{d}$. As in the proof
of Lemma~\ref{t1sg}, the key lies in defining  an appropriate rank
function $\rho$ on~$X$. Let $\Sigma_0$ be the set of vertices $x\in X$
that have finite degree in~$G$. For $\xi > 0$, let $\Sigma_{\xi}$ be
the set of vertices $x\in X$ such that any family of disjoint paths
starting at $x$ and ending in a vertex $y\not\in \bigcup_{\zeta\in\xi}
\Sigma_{\zeta}$ is finite. Again, we have $\Sigma_{\zeta}\subset
\Sigma_{\xi}$ for $\zeta < \xi$. Finally, for each $x\in X$, define
$\rho(x)$ to be the least $\xi$ such that $x\in\Sigma_{\xi}$; if no
such $\xi$ exists, let $\rho(x)$ remain undefined.

As in the proof of Lemma~\ref{t1sg}, we may imitate the standard
construction of~$T_\omega$ to show that if there exists an $x\in X$
such that $\rho(x)$ has remained undefined, then $G$ contains a
subdivision of~$T_\omega$.

We shall therefore assume that $\rho(x)$ is defined for all $x\in X$,
and show that $G$ is not dominating. Let $L\colon\,X\rightarrow \omega$
be any labelling. We shall find a function $H\colon\,\omega \to\omega$
which is not dominated by $L\circ R$ for any ray $R$ in~$G$.

Consider a vertex $x\in X$, and let $Y = \{ y : \rho(y)\geq\rho(x),\,
y\ne x \}$. Consider an a
rbitrary set $\cal P$ of disjoint paths
starting at $x$ and ending in a vertex of~$Y$. By the definition
of~$\rho$, any such set must be finite. As is easy to see, this
implies that there is in fact a common finite bound on the
cardinalities of all such sets~$\cal P$. Then $x$ must be separated
from~$Y$ by some finite set $Y_x\subset X\setminus\{x\}$ --- this
means that  every path from $x$ to a vertex of~$Y$ meets~$Y_x$ ---
because $Y_x$ can be chosen to be a maximal family of disjoint paths
starting at $x$ and ending at a vertex of $Y$.

For each $x\in X$, let $\{x\} = T_x^0 \subset T_x^1\subset T_x^2
\subset\ldots$ be an infinite sequence of finite subsets of~$X$, chosen
so that for every $i$ and $z\in T_x^i$ we have $Y_z\subset T_x^{i+1}$.
It  is then possible to define a function $Q_x\colon\,\omega\rightarrow
\omega$ such that $Q_x(m) \geq L(y)$ for every $m\in\omega$ and every
$y\in T_x^m$.  From our hypothesis that $\card{X} < \frak{d}$ it
follows that there exists a function
$H\colon\,\omega\rightarrow\omega$ which is not dominated by any of
the functions~$Q_x$; clearly, we may choose $H$ to be increasing.

Now let $R$ be any ray in~$G$; we prove that $H$ is not dominated by
$L\circ R$. By Lemma~\ref{tail}, there is some $K\in\omega$ such that
for each $i\geq K$ there is a $k > i$ with $\rho(R(i))\leq
\rho(R(k))$. Let
 $$M = \{ m\in\omega : H(m) > Q_{R(K)}(m)\}.$$ $M$~is infinite, since
$H\not\led  Q_{R(K)}$. We show that for each $m\in M$ with $m\geq K$
there is some $j\geq m$ such that $Q_{R(K)}(m) \geq L(R(j))$. Since
$H$ is increasing, this will imply that
 $$H(j)\geq H(m) > Q_{R(K)}(m) \geq L(R(j))$$
for all these infinitely many~$j$, giving $H\not\led  L\circ R$ as
desired.

It suffices to prove that  for each $m \geq K$ there is  some $j \geq
m$ such that $R(j) \in T_{R(K)}^{m-K}$ ($\subset T_{R(K)}^m$), because
then $Q_{R(K)}(m) \geq L(R(j))$ by definition. This fact can be proved
by induction on $m$. If $m = K$, let $j = K$; then $\{R(j)\} =
\{R(K)\} = T_{R(K)}^0 = T_{R(K)}^{m-K}$ as desired. If $m > K$, use
the induction hypothesis to find an $i \geq m-1$ such that $R(i) \in
T_{R(K)}^{m-1-K}$, and choose $k > i$ so that $\rho(R(i)) \leq
\rho(R(k))$. (Such $k$ exists by $m-1\ge K$ and the choice of~$K$.)
Then $Y_{R(i)}$ separates $R(i)$ from~$R(k)$, so there is a $j$ with
$i < j \leq k$ such that $R(j)\in Y_{R(i)}$. Then $R(j)\in
Y_{R(i)}\subset T_{R(K)}^{m-K}$ (by $R(i)\in T_{R(K)}^{m-1-K}$ and the
definition of~$T_{R(K)}^{m-K}$) and $j\geq i+1\geq m$, so $j$ is as
desired.
 \stopproof

Let us say that a function $f\colon\, \omega\to\omega$ {\em tends to
infinity} if $f^{-1} (n)$ is finite for every $n\in\omega$.

\begin{lemma}
    If $G$ is a graph on $X$, \label{submodel} and if $Y\subset X$
and $L\colon\, X\rightarrow\omega$, then there is a set $Z$ with
$Y\subset Z\subset X$ and $\card{Y} =\card{Z}$ which has the
following property: for any ray $R$ in $G$ with infinitely many
vertices in~$Z$ and $L\circ R$ tending to infinity, there is a
ray $R'$ in~$G\cap [Z]^2$ such that $L\circ R' = L\circ R$.
 \end{lemma}
 \proof
 The lemma is trivial when $Y$ is finite, so we assume that $Y$
is infinite. Beginning with $Z_0 = Y$, let us define an infinite
increasing sequence $Z_0\subset Z_1\subset Z_2\subset\dots$ of
subsets of~$X$, as follows. Suppose $Z_n$ has already been
defined. To obtain $Z_{n+1}$ from~$Z_n$, consider first every
vertex $y\in Z_n$. Let ${\cal P}$ be a maximal set of (finite)
paths in $G$ ending in $y$ and having no other vertices
in~$Z_n$ such that $L\circ P \ne L\circ P'$ for distinct $P,P'
\in {\cal P}$. (This implies that ${\cal P}$ is countable.) For
each $P\in {\cal P}$, check whether $G\cap [Z_n]^2$ contains an
infinite set of disjoint paths ending in~$y$ such that every
path $P'$ in this set satisfies $L\circ P' = L\circ P$; if there
is no such set then add the vertices of $P$ to~$Z_n$. Similarly,
consider every pair $\{x,y\} \in [Z_n]^2$. Now let ${\cal P}$ be
a maximal set of $x$--$y$ paths in $G$ that hav
e no other
vertices in~$Z_n$, and such that $L\circ P \ne L\circ P'$ for
distinct $P,P' \in {\cal P}$. For each $P\in {\cal P}$, check
whether $G\cap [Z_n]^2$ contains an infinite set of disjoint
paths from $x$ to~$y$ such that every path $P'$ in this set
satisfies $L\circ P' = L\circ P$; if there is no such set then
add the vertices of $P$ to~$Z_n$.

Note that, since $Y = Z_0$ was assumed to be infinite, we have
$\card{Z_n} = \card{Z_{n+1}}$ for each~$n$. Therefore $Z =
\bigcup_{n\in\omega} Z_n$ satisfies $\card{Y} = \card{Z}$ as
required. Moreover, $Z$~has the following two properties.
Whenever $y\in Z$ and $P$ is a path of length $>1$ in $G$ that
ends in $y$ but has no other vertices in~$Z$, there is an
infinite set of disjoint paths ending in~$y$ such that every
path $P'$ in this set has all its vertices in $Z$ and satisfies
$L\circ P' = L\circ P$. Similarly, whenever $x,y\in Z$ are
joined in $G$ by a path $P$ of length $>2$ whose only vertices
in $Z$ are $x$ and~$y$, there are infinitely many disjoint paths
$P'$ from $x$ to $y$ whose vertices are all in~$Z$ and which
satisfy $L\circ P' = L\circ P$.

Now let $R$ be any ray in $G$ with infinitely many vertices
in~$Z$ and $L\circ R$ tending to infinity. If all the vertices of
$R$ are in~$Z$, we set $R' = R$. Otherwise there is a (finite or
infinite) sequence $m_0 \le n_1 < m_1 \le n_2 < m_2 \le\dots$ of
integers such that the vertices of $R$ outside $Z$ are precisely
the vertices of the form $R(k)$ with $k < m_0$ or $n_i < k < m_i$
for some~$i$. We shall obtain $R'$ from $R$ by replacing its
initial segment $P_0 = R\upharpoonright m_0$ and, for $i>0$, its
subpaths $P_i$ from $x_i = R(n_i)$ to $y_i = R(m_i)$ with paths
on $Z$ that carry the same labelling.

For each $i = 0, 1, \dots$ in turn, let us find a path $Q_i$
in $G\cap [Z]^2$ from $x_i$ to~$y_i$ (or, in the case of $i=0$,
just ending in~$y_0 = R(m_0)$) so that $L\circ Q_i = L\circ P_i$.
If $P_i$ has no vertices outside~$Z$, we let $Q_i = P_i$.
Otherwise, by the construction of~$Z$, there is an infinite set
${\cal Q}_i$ of disjoint paths that qualify for selection
as~$Q_i$. Now ${\cal Q}_i$~has an infinite subset ${\cal Q}'_i$
of paths all avoiding the paths $Q_j$ chosen earlier (except
that we might have $x_i = y_j$ if $i = j+1$). Since the paths in
${\cal Q}'_i$ all carry the same labelling, they only use
finitely many labels. Since, by assumption, $L\circ R$ tends to
infinity, $R$~has only finitely many vertices carrying any of
these labels. Since ${\cal Q}'_i$ is an infinite set of disjoint
paths from $x_i$ to~$y_i$ (or ending at~$y_0$, respectively), we
may therefore choose $Q_i$ from ${\cal Q}'_i$ so that $Q_i$ has
no other vertices on~$R$.

Let $R'$ be obtained from $R$ by replacing each $P_i$ with the
corresponding $Q_i$ as defined above. Then $R'$ is a ray in
$G\cap [Z]^2$, and $L\circ R' = L\circ R$ as required.
 \stopproof

{\bf Proof of Theorem \ref{main}:} The sufficiency of the three
conditions has been established in Propositions \ref{type1},
\ref{type2} and \ref{rays}, respectively. To prove the
necessity, let $G$ be a dominating graph on a set $X$ and
suppose that this is witnessed by the function $L\colon
X\rightarrow \omega$. Let ${\cal R}$ be a maximal collection of
disjoint rays in $G$. If $\card{{\cal R}} \geq \frak{d}$ then
there is nothing to do. If not, it follows from Lemma
\ref{submodel} that there is some $Y\subset X$ such that
 \begin{itemize}
 \item $\card{Y} =\card{{\cal R}} < \frak{d}$;
 \item if $R\in{\cal R}$ then $R\subset Y$;
 \item for any ray $R$ in $G$ with infinitely many vertices
in~$Y$ and $L\circ R$ tending to infinity, there is a ray $R'$
in~$G\cap [Y]^2$ such that $L\circ R' = L\circ R$.
 \end{itemize} Let us show that $G\cap [Y]^2$ is a dominating
graph on~$Y$, and that this is witnessed by the
labelling~$L\upharpoonright Y$. Let $f\colon\, \omega\to\omega$
be given, without loss of generality increasing. Since $G$ is
dominating, it contains a ray $R$ such that $f\led L\circ R$.
Since $f$
 is increasing, $L\circ R$ tends to infinity. Moreover,
$R$~has infinitely many vertices in~$Y$, by the maximality
of~${\cal R}$. Therefore, $G\cap [Y]^2$ has a ray $R'$ such that
$L\circ R' = L\circ R$ and hence $f \led L\circ R'$.

Let $\cal T$ be a maximal collection of disjoint subdivisions
of~$T_\omega$ contained in $G\cap [Y]^2$. If $\card{{\cal T}}
\geq \frak{b}$ then there is nothing to do. If not, it follows
from Lemma~\ref{submodel} that there is some $Z\subset Y$ such
that
 \begin{itemize}
 \item $\card{Z} =\card{{\cal T}} < \frak{b}$;
 \item if $T\in{\cal T}$ then the vertices of $T$ are all in
$Z$;
 \item for any ray $R$ on $Y$ with infinitely many vertices
in~$Z$ and $L\circ R$ tending to infinity, there is a ray $R'$
in~$G\cap [Z]^2$ such that $L\circ R' = L\circ R$.
 \end{itemize} If  $G\cap [Z]^2$ contains a uniform subdivision
of~$T_\omega$, we are done; we therefore assume that it does
not. Then, by Lemma \ref{t1sg}, $G\cap [Z]^2$ is not dominating.
We show that now $G\cap [Y\setminus Z]^2$ must be a dominating
graph on~$Y\setminus Z$. Since $\card{Y\setminus Z} < \frak d$
and $G\cap [Y\setminus Z]^2$ contains no subdivision
of~$T_\omega$ (by the maximality of~$\cal T$), this will
contradict Lemma~\ref{t2sg}.

Let $H\colon\,\omega\to\omega$ be a function witnessing (with
respect to~$L$) that $G\cap [Z]^2$ is not dominating. In order
to show that $G\cap [Y\setminus Z]^2$ is dominating, let
$I\colon\,\omega \to \omega$ be given; we shall find a ray on
$Y\setminus Z$ whose sequence of labels dominates~$I$. Let
$J\colon\, \omega\to\omega$ be increasing and such that
$J(n)\ge\max\{H(n), I(n)\}$ for every~$n$. Recall that $G\cap
[Y]^2$ with $L$ was found to be dominating; choose a ray $R$ on
$Y$ so that $J\led  L\circ R$. As $J$ is increasing, $L\circ R$
tends to infinity. Since~$H$, and hence also~$J$, witnesses that
$G\cap [Z]^2$ is not dominating, the definition of $Z$ implies
that $R$ meets $Z$ in only finitely many vertices. Let $R'$ be a
subray of $R$ whose vertices are all in~$Y\setminus Z$; since
$J\led L\circ R$ and $J$ is increasing, we have $I\led  J\led 
L\circ R'$ as desired.
 \stopproof

\bigbreak\centerline{\bf ACKNOWLEDGEMENT}
 \medskip\noindent During their work on this paper, the authors
benefited from the continual interest of Paul Erd\H os and Alan
Mekler.

\bigbreak\bigskip\centerline{\bf REFERENCES} \medbreak

\REF 1 R.~Diestel and I.~Leader, A proof of the bounded graph
conjecture, submitted.

\REF 2 R.~Halin, Some problems and results in infinite graphs,
{\em Annals of Discrete Math.}~41 (1989), 195--210.

\REF 3 J.~van Mill and G.M.~Reed (editors), {\em Open Problems
in Topology}, North-Holland, Amsterdam, 1991.

\REF 4 F. Rothberger, Sur les Familles Ind\'{e}nombrables des
Suites de Nombres Naturels et les Probl\`{e}mes Concernants la
Propri\'{e}t\'{e} C, {\em Proceedings of the Cambridge
Philosophical Society} 27 (1941), 8--26.

\REF 5 F. Rothberger, Une Remarque Concernante l'Hypoth\`{e}se du
Continu, {\em Fundamenta Mathematicae} 31 (1938), 224-226.

\end{document}